\documentclass[12pt]{amsart}

\usepackage[latin1]{inputenc}
\usepackage{mathptmx}
\usepackage{amscd}

\newtheorem{theorem}{Theorem}[section]
\newtheorem{lemma}[theorem]{Lemma}
\newtheorem{corollary}[theorem]{Corollary}
\newtheorem{proposition}[theorem]{Proposition}

\theoremstyle{definition}
\newtheorem{remark}[theorem]{Remark}
\newtheorem{definition}[theorem]{Definition}
\newtheorem{example}[theorem]{Example}

\theoremstyle{plain}
\newtheorem*{thm*}{Theorem}

\newcommand{\mm}{\mathfrak m}

\newcommand{\pp}{\mathfrak p}
\newcommand{\qq}{\mathfrak q}
\def\NZQ{\Bbb}

\def\ZZ{{\NZQ Z}}
\def\RR{{\NZQ R}}

\def\NN{{\NZQ N}}

\newcommand{\Fcc}{\mathcal{F}}

\let\Bbb=\mathbb
\let\phi=\varphi

\def\relint{\operatorname{relint}}

\def\Hom{\operatorname{Hom}}
\def\Ext{\operatorname{Ext}}

\def\Ann{\operatorname{Ann}}

\def\Min{\operatorname{Min}}

\def\cha{\operatorname{char}}

\let\oldbigwedge\bigwedge
\def\BIGwedge{{\textstyle\oldbigwedge}}
\def\medwedge{{\scriptstyle\oldbigwedge}}
\def\bigwedge{\mathchoice{\BIGwedge}{\BIGwedge}{\medwedge}{}}

\let\hat=\widehat

\let\iso=\cong

\let\epsilon=\varepsilon
\let\tilde=\widetilde

\DeclareMathOperator{\cn}{cn}

\DeclareMathOperator{\Supp}{Supp}
\DeclareMathOperator{\fan}{Fan}
\DeclareMathOperator{\str}{star}
\DeclareMathOperator{\lk}{link}

\DeclareMathOperator{\lpnt}{\hbox{\large\bf.}}
\DeclareMathOperator{\pnt}{\raise 0.5mm \hbox{\large\bf.}}

\title{On canonical modules of toric face rings}

\author{Bogdan Ichim}
\address{Universit\"at Osnabr\"uck, FB Mathematik/Informatik, 49069
Osnabr\"uck, Germany;
\phantom{iii} Institute of Mathematics, C.P. 1-764, 70700 Bucharest, Romania}
\email{bogdan.ichim@math.uos.de, bogdan.ichim@imar.ro}

\author{Tim R\"omer}
\address{Universit\"at Osnabr\"uck, FB Mathematik/Informatik, 49069 Osnabr\"uck, Germany}
\email{troemer@uos.de}

\begin{document}
\begin{abstract}
Generalizing the 
concepts of Stanley--Reisner and affine monoid algebras,
one can associate to
a rational pointed fan $\Sigma$ in $\RR^d$
the $\ZZ^d$-graded toric face ring $K[\Sigma]$.
Assuming that $K[\Sigma]$ is Cohen--Macaulay,
the main result of this paper is
to characterize the situation when its canonical module
is isomorphic to a $\ZZ^d$-graded ideal of $K[\Sigma]$.
From this result several algebraic and combinatorial consequences are deduced
in the situations where $\Sigma$ may be related to a manifold with non-empty boundary,
or $\Sigma$ is a shellable fan.
\end{abstract}
\maketitle

%
%
%
\section{Introduction}
\label{intro}
Let $\Sigma$ be a rational pointed fan in $\RR^d$,
i.e.\ $\Sigma$ is a finite collection of rational pointed cones in $\RR^d$
such that for $C' \subseteq C$ with $C\in \Sigma$
we have that $C'$ is a face of $C$ if and only if $C' \in \Sigma$,
and if $C,C' \in \Sigma$, then $C\cap C'$ is a common face of $C$ and $C'$.
Stanley constructed in \cite{ST87}
the {\em toric face ring $K[\Sigma]$} over a field $K$ associated to $\Sigma$  as follows.

As a $K$-vectors space $K[\Sigma]=\bigoplus_{a \in \bigcup_{C\in \Sigma} C \cap \ZZ^d} Kx^a$.
The multiplication in $K[\Sigma]$ is defined by:
$$
x^a \cdot x^b =
\begin{cases}
x^{a+b} & \text{if $a$ and $b$ are elements of a common face $C\in \Sigma$,}\\
0 & \text{otherwise.}
\end{cases}
$$
We see that $K[\Sigma]$ is naturally a $\ZZ^d$-graded $K$-algebra.
This class of rings generalizes the concepts of Stanley--Reisner rings
associated to simplicial complexes and affine monoid algebras
associated to affine monoids in some $\ZZ^d$.
(See Bruns--Herzog \cite{BH} or Stanley \cite{ST05} for
a detailed discussion on these two special cases.)

Motivated by known
results from the two special cases, the authors started to study
toric face rings systematically to
obtain their algebraic properties in \cite{IR07} .
(See also  \cite{BBR05}, \cite{BRRO04}, \cite{BRRO05}, \cite{BRGUBbook} and \cite{YA06}  for related results.)

Two of the main results in \cite{IR07} classify completely
in which cases $K[\Sigma]$ is Cohen--Macaulay, respectively Gorenstein,
in combinatorial terms of the rational pointed fan $\Sigma$.
It is well-known that the Gorenstein property is  equivalent to the fact that $K[\Sigma]$ is Cohen--Macaulay and
some principal ideal $K[\Sigma](-\sigma)$ is isomorphic to $\omega_{K[\Sigma]}$ as a $\ZZ^d$-graded $K[\Sigma]$-module
for some $\sigma\in \ZZ^d$,
where $\omega_{K[\Sigma]}$ denotes the  $\ZZ^d$-graded canonical module of $K[\Sigma]$.
Assuming only that $K[\Sigma]$ is Cohen--Macaulay,
it is a natural question in which cases
$\omega_{K[\Sigma]}$ is isomorphic to a $\ZZ^d$-graded ideal of $K[\Sigma]$.
After introducing definitions and results related to toric face rings in Section \ref{toricfacerings},
we give a complete characterization of this situation in Section \ref{canonicalmodule}.

Intersecting each cone of $\Sigma$ with the unit sphere of $\RR^d$
one obtains a regular cell complex $\Gamma_\Sigma$.
In Section \ref{topology} we study algebraic properties
of $K[\Sigma]$ in terms of topological properties of $\Gamma_\Sigma$.

In particular, in Section \ref{topology}, we get as an application of our main result of Section
\ref{canonicalmodule} the following theorem which generalizes
a nice result of Hochster for  Stanley--Reisner rings.

\begin{theorem}
Let $\Sigma$ be a rational pointed fan in $\RR^d$ such that $K[\Sigma]$ is a Cohen--Macaulay ring.
Assume that $X=|\Gamma_\Sigma|$ is a manifold with a non-empty boundary $\partial X$.
Further let $\Sigma'$ be the subfan of $\Sigma$
such that $\partial X=|\Gamma_{\Sigma'}|$.
Then the following conditions are equivalent:
\begin{enumerate}
\item
$\omega_{K[\Sigma]}$
is isomorphic as a $\ZZ^d$-module to the kernel of the natural surjective
homomorphism $K[\Sigma] \to K[\Sigma']$;
\item
$K[\Sigma']$ is Gorenstein and
$\omega_{K[\Sigma']}\cong K[\Sigma']$ as $\ZZ^d$-modules.
\end{enumerate}
\end{theorem}
For the definition of an {\em Euler fan} we refer to  Section \ref{toricfacerings}.
Condition (ii) of the  theorem is equivalent to the fact that
$K[\Sigma']$ is Cohen--Macaulay and $\Sigma'$ is an Euler fan.

Recall
that a {\em shelling} of $\Sigma$
is a linear ordering $C_1,\dots, C_s$
of the facets of $\Sigma$ such that either $\dim \Sigma =0$,
or the following two conditions are satisfied:
\begin{enumerate}
\item
$\fan(\partial C_1)$ has a shelling.
\item
For $1<j\leq s$
there exists a shelling
$D_{1},\dots, D_{t_j}$ of the (pure) fan $\fan(\partial C_j)$
such that
$
\emptyset
\neq
\bigcup_{i=1}^{j-1} \fan(C_i) \cap \fan(C_j)
=
\bigcup_{l=1}^{r_j} \fan(D_l)$
for some $1\leq r_j \leq t_j$.
\end{enumerate}
$\Sigma$ is called {\em shellable} if it has a shelling.
In the last Section of this paper
we study algebraic properties
of the toric face rings
$K[\bigcup_{i=1}^{j-1} \fan(C_i) \cap \fan(C_j)]$.
We can prove that they are all Cohen--Macaulay
and using  the  result of Section \ref{topology} mentioned above,
we show that the
canonical modules
$\omega_{K[\bigcup_{i=1}^{j-1} \fan(C_i) \cap \fan(C_j)]}$
are isomorphic to $\ZZ^d$-graded ideals of
$K[\bigcup_{i=1}^{j-1} \fan(C_i) \cap \fan(C_j)]$ for all $1<j\le s$. Also remark 
that $K[\bigcup_{l=r_j+1}^{t_j} \fan(D_l)]=0$ or $\omega_{K[\bigcup_{l=r_j+1}^{t_j} \fan(D_l)]}$ is 
isomorphic to a $\ZZ^d$-graded ideal of $K[\bigcup_{l=r_j+1}^{t_j} \fan(D_l)]$
for all $1<j\le s$.
Moreover, this result holds also under the weaker assumption
``semishellable'' as defined in Section \ref{shellfan}.

For a finitely generated  $R$-module $M$ over a commutative ring $R$ we call a finite filtration
$$
0=M_0 \subset M_1 \subset \dots \subset M_r=M
$$
of submodules of $M$ a {\em prime filtration},
if for every $1\le i\le r$ there is an isomorphism $M_i/M_{i-1}\iso R/P_i$ for some prime ideal $P_i$ of $R$.
Such a prime filtration
is called {\em clean}
if its set of corresponding prime ideals is equal to $\Min(\Supp M)$.
Then we call the module $M$ {\em clean}, if $M$ admits a clean filtration.
For a Stanley--Reisner ring $K[\Delta]$ corresponding to a simplicial complex
$\Delta$, Dress \cite{D} proved that the shellability of $\Delta$
is equivalent to the fact that $K[\Delta]$ is a {\em clean} ring.
This result does not hold for toric face rings.
In Theorem \ref{clean} we prove that $K[\Sigma]$ is clean if and only if
$\Sigma$ is shellable and certain combinatorial conditions are satisfied.
If $K[\Sigma]$ is a clean ring, then there is a shelling $C_1,\dots, C_s$ of $\Sigma$ such that (with the above notation)  
$K[\bigcup_{l=r_j+1}^{t_j} \fan(D_l)]=0$ or $K[\bigcup_{l=r_j+1}^{t_j} \fan(D_l)]$ is Gorenstein 
for all $1<j\le s$. That is,  $K[\bigcup_{l=r_j+1}^{t_j} \fan(D_l)]=0$ or $\omega_{K[\bigcup_{l=r_j+1}^{t_j} \fan(D_l)]}$ is a principal ideal for all $1<j\le s$. This provides a way to distinguish between shellable and clean fans.

We are grateful to Prof.\ W. Bruns
for inspiring discussions on the subject of the paper.

%
%
%
\section{Toric face rings}
\label{toricfacerings}
In this section we fix some notation and recall important results needed
in the rest of the paper.
Let $\Sigma$ be a rational pointed fan in $\RR^d$.
Let $K$ be a field.
We consider
the toric face ring $K[\Sigma]$ of $\Sigma$ over $K$ as defined
in Section \ref{intro}.
For $C \in \Sigma$ we define the $\ZZ^d$-graded ideal $\pp_C=(x^a \in K[\Sigma] : a \not\in C)$.
Is is easily verified that $K[\Sigma]/\pp_C \cong K[\fan(C)]$ where
$\fan(C)$ is the fan consisting of all cones $D \in \Sigma$ such that $D \subseteq C$.
In particular, $\pp_C$ is a prime ideal.
More generally, for a subfan $\Sigma' \subseteq \Sigma$ we
set $\qq_{\Sigma'}=\qq_{\Sigma'}^{\Sigma}=(x^a \in K[\Sigma] : a \not\in |\Sigma'|)$.
As was observed in \cite{IR07},
all $\ZZ^d$-graded prime and radical ideals of $K[\Sigma]$ are described as follows.
\begin{lemma}
\label{primeradical_ideals}
Let $\Sigma$ be a rational pointed fan  in $\RR^d$.
\begin{enumerate}
\item
The assignment $C\mapsto \pp_C$
is a bijection between
the set of non-empty cones in $\Sigma$ and the set of $\ZZ^d$-graded prime ideals of $K[\Sigma]$.
\item
The assignment $\Sigma'\mapsto \qq_{\Sigma'}$
is a bijection between
the set of non-empty subfans of $\Sigma$
and the set of $\ZZ^d$-graded radical ideals of $K[\Sigma]$.
\end{enumerate}
In particular, $\mm=\pp_0$ is the unique $\ZZ^d$-graded maximal ideal of $K[\Sigma]$,
which is also maximal in the usual sense.
\end{lemma}

For $C \in \Sigma$ let $\str_\Sigma(C)=\{D \in \Sigma : C \subseteq D \}$
be the {\em star of $C$}. In particular, $\str_\Sigma(0)=\Sigma$.
Moreover, we set $\Sigma(C)=\Sigma\setminus \str_\Sigma(C)$.
Observe that this a subfan of $\Sigma$.

We define a complex

$$
\CD
\mathcal{C}^{\lpnt}_{\str_\Sigma(C)}
\colon
0@>>>\mathcal{C}^{\dim C-1}@>\partial>>\mathcal{C}^{\dim C}@>>>\cdots @>\partial>>\mathcal{C}^{\dim\Sigma-1}@>>>0
\endCD
$$
where
$$
\mathcal{C}^{i}=\bigoplus_{D\in \str_\Sigma(C),\ \dim D=i+1}K D
$$
for $i=\dim C-1,\ldots,\dim\Sigma-1$,
and  the differential is induced by an incidence function $\varepsilon$ on $\Sigma$.
(See \cite[Section 4]{IR07} for details. There $\mathcal{C}^{\lpnt}_{\str_\Sigma(C)}$ was denoted by
$\mathcal{C}^{\lpnt}(\Gamma_{\str_\Sigma(C)})$.)
Note that the (co-)homology
$\Tilde H^{i}(\mathcal{C}^{\lpnt}_{\str_\Sigma(C)})$
of the complex $\mathcal{C}^{\lpnt}_{\str_\Sigma(C)}$
computes
up to a shift the reduced cohomology groups
$\Tilde H^{i}(\Delta(\str_\Sigma(C)))$ with coefficients in $K$
of the order complex
$\Delta(\str_\Sigma(C))$ of the poset $\str_\Sigma(C)-\{C\}$.
Indeed, in \cite[Lemma 4.6]{IR07} it was shown that
$\Tilde H^{i}(\mathcal{C}^{\lpnt}_{\str_\Sigma(C)}) =
\Tilde H^{i-\dim C}(\Delta(\str_\Sigma(C)))$
for all integers $i$.

Observe that for $a\in \relint(C)\cap \ZZ^d$ we have
$\str_\Sigma(C)=\{D \in \Sigma : a \in D\}$.
Hence it makes sense to define slightly more generally for $a\in \ZZ^d$ the {\em star of $a$} as
$\str_\Sigma(a)=\{D\in\Sigma:a\in D\}$. Analogously we define $\Sigma(a)$ and the complex
 $\mathcal{C}^{\lpnt}_{\str_\Sigma(a)}$.
By the local cohomology groups of $K[\Sigma]$
we always mean the local cohomology with respect to the
$\ZZ^d$-graded maximal  ideal $\mm$  of $K[\Sigma]$.
A Hochster type formula for the local cohomology of toric face rings is given by:
\begin{theorem}[Corollary 4.7. \cite{IR07}]
\label{hochster}
Let $\Sigma$ be a rational pointed fan  in $\RR^d$.
Then
\begin{eqnarray*}
H_\mm^i(K[\Sigma])
&\iso&
\bigoplus_{C\in \Sigma}
\bigoplus_{a \in -\relint(C)}
\Tilde H^{i-1}(\mathcal{C}^{\lpnt}_{\str_\Sigma(C)}) \otimes_K K(-a)\\
\end{eqnarray*}
as $\ZZ^d$-graded $K$-modules (where
$\relint(C)$ denotes the relative interior of $C$ with respect to the subspace topology on the vector space generated by $C$).
\end{theorem}

See also \cite{BBR05} where such a formula was given for a more general class of rings.
Theorem \ref{hochster} yields an important Cohen--Macaulay criterion for fans
which generalizes the corresponding result for simplicial complexes given by Reisner in \cite{RE76}:

\begin{corollary}
\label{reisner}
 Let $\Sigma$ be a rational pointed fan  in $\RR^d$.
Then $K[\Sigma]$ is Cohen--Macaulay if and only if $\Tilde H^{i}(\mathcal{C}^{\lpnt}_{\str_\Sigma(C)})=0$
for all $C\in\Sigma$ and all $i<\dim \Sigma-1.$
\end{corollary}

Generalizing ideas of \cite{GR84} and \cite{GR842} we consider the complex
$$
\CD
D_{\pnt}(\Sigma): 0@>>> D_{\dim\Sigma} @>\partial>> D_{\dim\Sigma-1}@>>>\cdots @>>> D_{1} @>\partial>> D_0 @>>> 0,
\endCD
$$
 where we set
$$
D_i=\bigoplus_{C\in \Sigma,\ \dim C=i}K [C\cap \ZZ^d]
$$
for $i=0,\ldots,\dim\Sigma$,
and  the differential is the canonical projection multiplied with an incidence function $\varepsilon$ on $\Sigma$.
(We refer again to \cite[Section $5$]{IR07} for a comprehensive treatment.)

Assume that $K[\Sigma]$ is Cohen--Macaulay.
Since $K[\Sigma]$ is a factor ring of a polynomial ring,
the canonical module $\omega_{K[\Sigma]}$ exists.
The module $\omega_{K[\Sigma]}$ is described by:

\begin{theorem}[Corollary 5.2. \cite{IR07}]
\label{simplified_DC}
Let $\Sigma$ be a rational pointed fan in $\RR^d$ with $k = \dim \Sigma$
and assume that $K[\Sigma]$ is Cohen--Macaulay.
Then $\omega_{K[\Sigma]}$ is isomorphic to the kernel of
$D_{k} \overset{\partial}{\to}  D_{k-1}$
and
the complex of $\ZZ^d$-graded $K[\Sigma]$-modules
$$
\CD
 0@>>> \omega_{K[\Sigma]}@>>>D_{k} @>\partial>> D_{k-1} @>>>\cdots @>>> D_{1} @>\partial>> D_0 @>>> 0
\endCD
$$
is exact.
\end{theorem}

For a cone $C\in\Sigma$
let $f_i(\str_\Sigma(C))$ be the number of $i$-dimensional
cones in $\str_\Sigma(C)$.
We define
$$
\rho_\Sigma (C)
=
\sum_{i=0}^{\dim \Sigma} (-1)^{i+1} f_i(\str_\Sigma(C)).
$$
Observe that the definition of $\rho$ differs from the definition of
the number $\Tilde\chi_\Sigma (C)$ as considered in \cite{IR07}
by the sign $(-1)^{\dim C+1}$.
Hence all results cited from \cite{IR07} have to take into account this factor.
But we will see soon that one gets more compact formulas using
$\rho_\Sigma (C)$ as defined above.
An easy computations shows that
$$
\rho_\Sigma (C)
=
\sum_{i=\dim C-1}^{\dim \Sigma -1 } (-1)^{i} \dim_K (\mathcal{C}^{i}_{\str_\Sigma(C)})
=
\sum_{i=\dim C-1}^{\dim \Sigma -1 } (-1)^{i} \dim_K \Tilde H^{i}(\mathcal{C}^{\lpnt}_{\str_\Sigma(C)}),
$$
so $\rho_\Sigma (C)$ is nothing else that the {\it Euler characteristic} of the complex $\mathcal{C}^{\lpnt}_{\str_\Sigma(C)}$.

Gorenstein  toric face rings are described  by the next theorem.

\begin{theorem}[Theorem 6.2. \cite{IR07}]
\label{gorensteinsecond}
Let $\Sigma$ be a rational pointed fan in $\RR^d$.
Then  $K[\Sigma]$ is Gorenstein if and only if $K[\Sigma]$ is Cohen--Macaulay and there
exists $\sigma\in\big [\bigcap_{C \in \Sigma \text{ maximal}}C\big ]\cap\ZZ^d$ such that:
\begin{enumerate}
\item
$$
\bigcup_{C\in \str_\Sigma(\sigma)}\relint(C)\cap \ZZ^d
=
\sigma+\Supp(K[\Sigma],\ZZ^d).$$
\item
For all cones $C \in \Sigma$
we have
$$
\rho_\Sigma (C)=
\begin{cases}
(-1)^{\dim \Sigma - 1} &\text{ if } C \in \str_\Sigma(\sigma) ;\\
0 &\text{ if } C \in \Sigma(\sigma).
\end{cases}
$$
\end{enumerate}
\end{theorem}

We call
$\Sigma$ an {\em Euler fan} if $\Sigma$ is pure
(i.e.\ all facets of $\Sigma$ have the same dimension),
and $\rho_\Sigma (C)=(-1)^{\dim \Sigma - 1}$ for all $C \in \Sigma$.
We have the following characterizations of Cohen--Macaulay
and Euler fans:
\begin{proposition}[Corollary 6.5. \cite{IR07}]
\label{needit}
Let $\Sigma$ be a rational pointed fan in $\RR^d$.
The following statements are equivalent:
\begin{enumerate}
\item
$K[\Sigma]$ is Cohen--Macaulay and $\Sigma$ is an Euler fan;
\item
$K[\Sigma]$ is
Gorenstein and $\omega_{K[\Sigma]} \cong K[\Sigma]$ as $\ZZ^d$-graded modules;
\item
$K[\Sigma]$ is Gorenstein and $\rho_\Sigma (0)=(-1)^{\dim \Sigma-1}$;
\item
$K[\Sigma]$ is Cohen--Macaulay and $\Tilde H^{\dim \Sigma-1}(\mathcal{C}^{\lpnt}_{\str_\Sigma(C)})=K$
for all $C \in \Sigma$.
\end{enumerate}
\end{proposition}

%
%
%
\section{Embeddings of the canonical module
as an ideal of a toric face ring}
\label{canonicalmodule}

Let $K$ be a field, and $\Sigma$ be a rational pointed fan in $\RR^d$ such that $K[\Sigma]$ is a Cohen--Macaulay ring.
An {\em admissible $\ZZ$-grading} on $K[\Sigma]$ is a $\ZZ$-grading such that
the ring $K[\Sigma]=\bigoplus_{i\in \NN} K[\Sigma]_i$ is a finitely generated positively graded $K$-algebra with
$K=K[\Sigma]_0$ and all components $K[\Sigma]_i$ are direct sums of
finitely many $\ZZ^d$-graded components.
Observe that there exists fan's such that there exists no such gradings on $K[\Sigma]$. 
(See \cite[Example 2.7]{BRGU01}.)

In \cite[Corollary 5.3.]{IR07}  it was shown
that in case that there exists an admissible $\ZZ$-grading on $K[\Sigma]$,
then there exists an  $\ZZ$-graded embedding of the canonical module $\omega_{K[\Sigma]}$ into $K[\Sigma]$.
Then $\omega_{K[\Sigma]}$ can be identified with a $\ZZ$-graded ideal $I$ of $K[\Sigma]$.

In general,
if $\omega_{K[\Sigma]}$ can be identified with an ideal $I$ of $K[\Sigma]$
we cannot expect that $I$ is $\ZZ^d$-graded,
simply because it may happen that $\dim_K (\omega_{K[\Sigma]})_a>1$ for some $a\in \ZZ^d$.
See \cite[Section 5.7]{BH} for such an example for a Stanley--Reisner ring.
The case when $I$ is $\ZZ^d$-graded is however interesting and
it is a natural question to characterize this situation.

First we need a special construction.
We always embed $\RR^d$ naturally into $\RR^{d+1}$ onto the first $d$ coordinates
and identify the basis-vectors $e_1,\dots,e_d$ of the corresponding standard basis
of $\RR^d$ and $\RR^{d+1}$ respectively.
Let  $C$ be a rational pointed cone in $\RR^d$ and $w$ a rational point in $\RR^{d+1}\setminus\RR^d$
(e.g.\ $w=e_{d+1}$).
We denote by $w\ast C$ the smallest (pointed) cone containing both $w$ and $C$. If $\Sigma$ is  rational pointed fan in $\RR^d$,
we denote by $w\ast \Sigma$ the fan with facets
$\{w\ast C: C \text{ facet of  }\Sigma \}$. Clearly $w\ast \Sigma$ is a
$(\dim \Sigma +1)$-dimensional fan in $\RR^{d+1}$.
For a $\ZZ^d$-graded module $M$ we let
$\Supp(M,\ZZ^d)=\{a \in \ZZ^d:M_a\not= 0\}$ be the {\em support} of $M$ in $\ZZ^d$.
Now we prove the following theorem which generalizes a well-known result for Stanley--Reisner rings.

\begin{theorem}
\label{main}
Let $\Sigma$ be a rational pointed fan in $\RR^d$ such that $K[\Sigma]$ is a Cohen--Macaulay ring.
Then the following statements are equivalent:
\begin{enumerate}
\item
$\Sigma$ is not an Euler fan, and there exists an embedding $\omega_{K[\Sigma]}\to K[\Sigma]$ of $\ZZ^d$-graded $K[\Sigma]$-modules;
\item
There exists a $(\dim\Sigma-1)$-dimensional subfan
$\Sigma'$ of $\Sigma$ which is Euler and Cohen--Macaulay over $K$ such that for all
$C\in\Sigma$
$$
 \Tilde H^{\dim\Sigma-1}(\mathcal{C}^{\lpnt}_{\str_\Sigma(C)}) =
\begin{cases}
0 & \text{if }  C\in\Sigma',\\
K & \text{if }  C\not\in\Sigma'.
\end{cases}
$$
\end{enumerate}
In particular, if the image of the embedding $\omega_{K[\Sigma]}\to K[\Sigma]$ is a principal ideal, then $K[\Sigma]$ is Gorenstein.
\end{theorem}
\begin{proof}
Before proving the equivalence of (i) and (ii),
we observe the following.
For
a $\ZZ^d$-graded $K[\Sigma]$-module $M$ we
let $M^\vee=\Hom_K(M,K)$ be the $\ZZ^d$-graded $K$-dual of $M$ which is naturally
again a $\ZZ^d$-graded $K[\Sigma]$-module with homogeneous components $(M^\vee)_a=\Hom_K(M_{-a},K)$ for $a \in \ZZ^d$.
It follows from the local duality theorem for $\ZZ^d$-graded modules
(see \cite[Theorem 3.6.19.]{BH} for the $\ZZ$-graded case) that
$$
\omega_{K[\Sigma]}\iso (H_\mm^{\dim\Sigma}(K[\Sigma]))^\vee\iso
\bigoplus_{a \in \ZZ^d}\Hom_K(H_\mm^{\dim\Sigma}(K[\Sigma])_{-a},K).
$$
Using this fact and
Theorem \ref{hochster} we deduce that
there exists an isomorphism of
$\ZZ^d$-graded $K$-vector spaces:
\begin{equation}
\label{helphochster}
\omega_{K[\Sigma]}\iso\bigoplus_{C\in \Sigma}\bigoplus_{a \in \relint(C)\cap \ZZ^d}
\Tilde H^{\dim\Sigma-1}(\mathcal{C}^{\lpnt}_{\str_\Sigma(C)}) \otimes_K K(-a).
\end{equation}
(i)$\Rightarrow$(ii): Let $I$ be the $\ZZ^d$-graded ideal which is isomorphic to the image of $\omega_{K[\Sigma]}\to K[\Sigma]$.
Since we assume that $\Sigma$ is not an Euler fan, it follows from
Proposition \ref{needit} that $I\not=K[\Sigma]$.
It follows from (\ref{helphochster}) that for
$C\in \Sigma$ we have either
$$
\relint(C)\cap \ZZ^d \subseteq \Supp(I,\ZZ^d)
\text{ or }
\relint(C)\cap \ZZ^d \bigcap \Supp(I,\ZZ^d)=\emptyset.
$$
We define
\begin{eqnarray*}
\Sigma'
&=&
\{C\in \Sigma: \relint(C)\cap \ZZ^d \bigcap \Supp(I,\ZZ^d)=\emptyset\}\\
&=&
\{C\in \Sigma:\Tilde H^{\dim\Sigma-1}(\mathcal{C}^{\lpnt}_{\str_\Sigma(C)})=0\}
\end{eqnarray*}
and claim that $\Sigma'$ is a subfan of $\Sigma$. Let $C\in\Sigma'$ and $C'$ be a face of $C$.
Assume that $C'\not\in\Sigma'$ and take $a,a'\in\ZZ^d$ such that $a\in \relint(C)$ and $a'\in \relint(C')$.
Then $x^a\not\in I$ and $x^{a'}\in I$. Since $a+a'\in \relint(C)$,
it follows that $x^a x^{a'}=x^{a+a'}\not\in I$, a contradiction
since $I$ is an ideal.
Thus we see that indeed  $\Sigma'$ is a subfan of $\Sigma$.
Moreover, $I=\qq_{\Sigma'}$ and by definition we have
for all
$C\in\Sigma'$
$$
 \Tilde H^{\dim\Sigma-1}(\mathcal{C}^{\lpnt}_{\str_\Sigma(C)}) =
\begin{cases}
0 & \text{if }  C\in\Sigma',\\
K & \text{if }  C\not\in\Sigma'.
\end{cases}
$$
Since $I\not=K[\Sigma]$ we deduce $\dim\Sigma>0$.
If $C\in\Sigma$ is a facet, then
$\Tilde H^{\dim\Sigma-1}(\mathcal{C}^{\lpnt}_{\str_\Sigma(C)})\not=0$. For each facet $C$ of $\Sigma$ we choose
$a_C\in \relint(C)\cap \ZZ^d$.
It follows from
the definition of the ring structure of $K[\Sigma]$
that
$$
\sum_{C\in \Sigma\text{ facet}}x^{a_C}\in I
$$
is a $K[\Sigma]$-regular element.
$I\cong \omega_{K[\Sigma]}$ is a canonical module for $K[\Sigma]$
in the category of  $K[\Sigma]$-modules.
An ideal which contains a regular element has a rank,
and we may use the $\ZZ^d$-graded analogue of \cite[Proposition 3.3.18 (b)]{BH}
to see that
$K[\Sigma']=K[\Sigma]/I$ is Gorenstein (in particular Cohen--Macaulay) of dimension $\dim\Sigma-1$.
Applying the functor
$\Hom_{K[\Sigma]}(\_, \omega_{K[\Sigma]})$
in the category of $\ZZ^d$-graded $K[\Sigma]$-modules to the exact sequence of $\ZZ^d$-graded $K[\Sigma]$-modules
$$
\CD
0@>>>\omega_{K[\Sigma]}@>>>K[\Sigma]@>>>K[\Sigma']@>>>0,
\endCD
$$
we obtain the exact sequence
\begin{align*}
0\longrightarrow\omega_{K[\Sigma]}\longrightarrow\Hom_{K[\Sigma]}(\omega_{K[\Sigma]},\omega_{K[\Sigma]})
&\longrightarrow\Ext^1_{K[\Sigma]}(K[\Sigma'],\omega_{K[\Sigma]})\\
&\longrightarrow\Ext^1_{K[\Sigma]}(K[\Sigma],\omega_{K[\Sigma]})\longrightarrow\cdots
\end{align*}
where we used the facts that
$\Hom_{K[\Sigma]}(K[\Sigma'],\omega_{K[\Sigma]})=0$ and
$\Hom_{K[\Sigma]}(K[\Sigma],\omega_{K[\Sigma]})=\omega_{K[\Sigma]}$.
Furthermore we have $\Hom_{K[\Sigma]}(\omega_{K[\Sigma]},\omega_{K[\Sigma]})\iso K[\Sigma]$
and
$\Ext^1_{K[\Sigma]}(K[\Sigma'],\omega_{K[\Sigma]})\iso
\omega_{K[\Sigma']}$
as $\ZZ^d$-graded $K[\Sigma]$-modules.
We also have $\Ext^1_{K[\Sigma]}(K[\Sigma],\omega_{K[\Sigma]})=0$.
(See \cite{BH}
for the $\ZZ$-graded cases of these results.)

So there exists a short exact sequence of $\ZZ^d$-graded $K[\Sigma]$-modules
$$
\CD
0@>>>\omega_{K[\Sigma]}@>>>K[\Sigma]@>>>\omega_{K[\Sigma']}@>>>0.
\endCD
$$
We conclude that
$K[\Sigma']\iso\omega_{K[\Sigma']}$ as $\ZZ^d$-graded modules and it follows now from
Proposition \ref{needit} that $\Sigma'$ is an Euler fan.

(ii)$\Rightarrow$(i):
$\Sigma'\not=\emptyset $ and $\Tilde H^{\dim \Sigma-1}(\mathcal{C}^{\lpnt}_{\str_\Sigma(C)})=0$
for all $C \in \Sigma'$. By Proposition \ref{needit} (iv), $\Sigma$ is not an Euler fan.

Now we show that $\qq^\Sigma_{\Sigma'}$ is isomorphic with the canonical module of $\Sigma$ which proves (i).
Let $w$ a rational point in $\RR^{d+1}\setminus\RR^d$ and
consider
the rational pointed fan $\Pi=(w\ast\Sigma') \cup\Sigma$ in $\RR^{d+1}$
with facets the facets of the fan $w\ast\Sigma'$ and the facets of the fan $\Sigma$.
(Here we consider $\Sigma$ as the fan in $\RR^{d+1}$ induced by the embedding $\RR^d\to \RR^{d+1}$.)

We claim that $\Pi$
is Euler and $K[\Pi]$ is Cohen--Macaulay of dimension $\dim\Sigma$.  For this we prove that for all
$C\in\Pi$
\begin{equation}
\label{eq}
\Tilde H^{i}(\mathcal{C}^{\lpnt}_{\str_\Pi(C)}) =
\begin{cases}
0 & \text{if }  i<\dim\Sigma-1,\\
K & \text{if }  i=\dim\Sigma-1.
\end{cases}
\end{equation}
Then it follows
from  Theorem \ref{reisner} that $K[\Pi]$ is Cohen--Macaulay of dimension $\dim\Sigma$
and using (iv) in Proposition \ref{needit} we can conclude that $\Pi$ is an Euler fan.

We consider the following cases:

(a) If $C\in\Sigma\setminus\Sigma'$, then $\str_\Pi(C)=\str_\Sigma(C)$ and
(\ref{eq}) follows from the assumption of (ii) and the
Cohen--Macaulayness of $K[\Sigma]$ (using also Corollary \ref{reisner}).

(b) If $w\in C$, then $C=w\ast F$
where $F\in\Sigma'$ and thus $\str_\Pi(C)$ is isomorphic as posets to $\str_{\Sigma'}(F)$.
Thus (\ref{eq}) follows since $\Sigma'$  is Euler and Cohen--Macaulay over $K$.
(One has to use again Corollary \ref{reisner} and Proposition \ref{needit} applied to $K[\Sigma']$.)

(c) If $C\in \Sigma'$,
then $\str_\Pi(C)=\{w\ast F:F\in\str_{\Sigma'}(C)\}\cup\str_\Sigma(C)$.
Moreover,
we have that $\{w\ast F:F\in\str_{\Sigma'}(C)\}\cap\str_\Sigma(C)=\emptyset$.
Set $\mathcal{C}^{\lpnt}_{\str_{\Sigma'}(C)}[-1]$ to be
the complex $\mathcal{C}^{\lpnt}_{\str_{\Sigma'}(C)}$ right shifted by one position
(that is $\mathcal{C}^{i}_{\str_{\Sigma'}(C)}[-1]=
\mathcal{C}^{i-1}_{\str_{\Sigma'}(C)}$). Then  we get an exact sequence of complexes
$$
\CD
0@>>>\mathcal{C}^{\lpnt}_{\str_\Sigma(C)}
@>>>\mathcal{C}^{\lpnt}_{\str_\Pi(C)}
@>>>\mathcal{C}^{\lpnt}_{\str_{\Sigma'}(C)}[-1]@>>>0.
\endCD
$$
This yields the long exact cohomology sequence
$$
\CD
\cdots
@>>>\Tilde H^{i}(\mathcal{C}^{\lpnt}_{\str_\Sigma(C)})
@>>>\Tilde H^{i}(\mathcal{C}^{\lpnt}_{\str_\Pi(C)})
@>>>\Tilde H^{i-1}(\mathcal{C}^{\lpnt}_{\str_{\Sigma'}(C)})
@>>>\cdots.
\endCD
$$
Since $K[\Sigma]$ is Cohen--Macaulay, it follows from
Corollary \ref{reisner}
that
$\Tilde H^{i}(\mathcal{C}^{\lpnt}_{\str_\Sigma(C)})=0$ for all $i\neq \dim \Sigma-1$.
But $C\in\Sigma'$ and hence also
$\Tilde H^{\dim \Sigma-1}(\mathcal{C}^{\lpnt}_{\str_\Sigma(C)})=0$.
We conclude that
$$
\Tilde H^{i}(\mathcal{C}^{\lpnt}_{\str_\Pi(C)})
\iso
\Tilde H^{i-1}(\mathcal{C}^{\lpnt}_{\str_{\Sigma'}(C)})
\quad\quad \text{for all } i.
$$
As $\Sigma'$  is Euler and $K[\Sigma']$ is Cohen--Macaulay,
we immediately see that (\ref{eq}) holds for $\Pi$.

We have shown that $K[\Pi]$ is Cohen--Macaulay and $\Pi$ is Euler. Thus it follows from Proposition \ref{needit}
that $K[\Pi]$ is Gorenstein and $\omega_{K[\Pi]} \cong K[\Pi]$ as $\ZZ^d$-graded modules.
Then the $\ZZ^d$-graded analogue if \cite[Proposition 3.6.12]{BH} implies the
$\ZZ^d$-graded
$K[\Sigma]$-modules isomorphisms
$$
\omega_{K[\Sigma]}\iso\Hom_{K[\Pi]}(K[\Sigma],K[\Pi])\iso\Ann_{K[\Pi]}\qq_\Sigma^\Pi.
$$
The monomials in the ideal $\Ann_{K[\Pi]}\qq_\Sigma^\Pi$ of $K[\Pi]$
are precisely those $x^a$ with
$a \in \Sigma\setminus \Sigma'$. Thus as a
$\ZZ^d$-graded $K[\Sigma]$-module it is isomorphic to $\qq_{\Sigma'}^\Sigma$.
Hence
$$
\omega_{K[\Sigma]} \iso \qq_{\Sigma'}^\Sigma
$$
as desired.
Finally, note that if $\qq_{\Sigma'}^\Sigma$ is principal,
it must be generated by a homogeneous non-zero divisor, because
$\omega_{K[\Sigma]}$ is a faithful $K[\Sigma]$-module.
So $\qq_{\Sigma'}^\Sigma \iso K[\Sigma](-a)$ for some $a \in \ZZ^d$
and thus $K[\Sigma]$ is Gorenstein.
\end{proof}

%
%
%
\section{On manifolds related to rational pointed fans}
\label{topology}
For further results on toric face rings we need
some facts related to a regular cell complex induced by a rational pointed fan $\Sigma$ in $\RR^d$.
For this we consider the intersection of $\Sigma$ with the unit sphere $\Bbb{S}^{d-1}$ and the set
$\Gamma_\Sigma=\big\{\relint(C)\cap \Bbb{S}^{d-1} : C\in\Sigma\big\}$.
For a cone $C\in \Sigma$ we denote by $e_C=\relint(C)\cap \Bbb{S}^{d-1}$ the corresponding element of $\Gamma_\Sigma$.
The elements $e_C$ are called {\em open cells}.
Set $\Gamma_\Sigma^i=\big\{e_C \in \Gamma_\Sigma:\bar e_C \text{ homeomorphic with }\Bbb{B}^i\big\}$
where $\Bbb{B}^i$ is the $i$-dimensional ball in $\RR^i$.
Now define $|\Gamma_\Sigma|=\bigcup_{C \in \Gamma} e_C$. Thus $|\Gamma_\Sigma|$ is just the intersection of
the underlying topological space $|\Sigma|=\bigcup_{C \in \Sigma} C$ of $\Sigma$ with the unit sphere $\Bbb{S}^{d-1}$.
Then $(|\Gamma_\Sigma|,\Gamma_\Sigma)$ is a {\em finite regular cell complex}.
The dimension of $\Gamma_\Sigma$ is given by
$\dim \Gamma_\Sigma=\max\{i:\Gamma_\Sigma^i\not=\emptyset\}= \dim \Sigma -1$.
An element $e_{C'}$ is called a {\em face} of $e_C$ if $e_{C'} \subset \bar e_C$, i.e.
$C'$ is a face of $C$.
Note that we consider $\emptyset=e_{\{0\}}$ also as a cell in $\Gamma_\Sigma$.

The following result is crucial.
First we need some further notations.
We already considered simplicial complexes $\Delta$
on a finite vertex set $V$. The elements $F \in \Delta$
are called {\em faces} of $\Delta$.
If $F \neq \emptyset$ we define the {\em dimension} of $F$ to be  $\dim F=\# F-1$. We set $\dim \emptyset =-1$ and
then the dimension $\dim \Delta$ of $\Delta$ is
the maximum of $\dim F$ for $F \in \Delta$.
For $F \in \Delta$ let $\lk_\Delta F$ be the set
$\{G \in \Delta: F\cup G \in \Delta,\ F\cap G =\emptyset\}$
and call it the {\em link} of $F$.
Simplicial complexes (and their faces)
have geometrically realizations which we usually denote by $|\Delta|$
and then $\Delta$ is said to be a triangulation of $|\Delta|$.
We recall one important construction. Let $P$ be a finite partially ordered set (poset for short) which has a unique minimal element $\hat 0$
(e.g.
the face poset of $\Sigma$,
the face poset of $\Gamma_\Sigma$ or $\str_\Sigma C$).
We denote by $\Delta(P)$ the {\em order complex}
of $P\setminus\{\hat 0\}$
which is the simplicial complex
where the faces are the chains of $P\setminus\{\hat 0\}$.
For further definitions and results on simplicial complexes  and posets see \cite{BH} and \cite{ST99} respectively.

Furthermore, for a topological space $X$, $p \in X$ and $K$ a field
we denote by $\tilde H_i(X)$ the reduced simplicial homology with coefficients in $K$, and by $H_i(X,X\setminus \{p\})$ the
local homology groups of $X$ at $p$ with coefficients in $K$.

\begin{lemma}
\label{topolinterpreter}
Let $\Sigma$ be a rational pointed fan in $\RR^d$,
$0\not= C \in \Sigma$ and $p \in e_C$.
Then
$$
H_i(|\Gamma_\Sigma|, |\Gamma_\Sigma| \setminus\{p\})
\cong\Tilde H^i(\mathcal{C}^{\lpnt}_{\str_\Sigma C})
\text{ for all } i.
$$
Moreover,
$$
\tilde H_i(|\Gamma_\Sigma|)
\cong\Tilde H^i(\mathcal{C}^{\lpnt}_{\str_\Sigma (0)})
\text{ for all } i.
$$
\end{lemma}
\begin{proof}
Let $\Fcc$ be a maximal chain
$C_1 \subset \dots \subset C_r$
in the order complex $\Delta(\Sigma)$ of the face poset of $\Sigma$  such that $C_r=C$.
Observe that then $\dim \Fcc= \dim C -1$
and
$\lk_{\Delta(\Sigma)} \Fcc = \Delta(\str C)$.
Note further that
$\Delta(\Sigma)$ is a triangulation of $\Gamma_\Sigma$
(see \cite[Proposition 4.7.8]{BJetal99}) and has therefore
the geometric realization $|\Delta(\Sigma)|=|\Gamma_\Sigma|$.
In particular, we may assume $p \in |\Delta(\Sigma)|$.
We compute
\begin{eqnarray*}
\Tilde H^i(\mathcal{C}^{\lpnt}_{\str_\Sigma C})
&\cong&
\Tilde H^{i-\dim C}(\Delta(\str_\Sigma C))\\
&\cong&
\Tilde H_{i-\dim C}(\Delta(\str_\Sigma C))\\
&=&
\Tilde H_{i-\dim C}(\lk_{\Delta(\Sigma)} \Fcc ) \\
&\cong&
H_{i}(|\Delta(\Sigma)|,|\Delta(\Sigma)|\setminus \{p\}) \\
&\cong&
H_{i}(|\Gamma_\Sigma|,|\Gamma_\Sigma|\setminus \{p\}).
\end{eqnarray*}
Here the first isomorphism was shown in \cite[Lemma 4.6]{IR07}.
The second isomorphism follows since we take (co-) homology with coefficients in $K$.
The third equality follows from the observations above.
The forth isomorphism is due to \cite[Lemma 5.4.5]{BH} . Here note that if $\str_\Sigma C = \{C\}$
(that is $\lk_{\Delta(\Sigma)} \Fcc =\emptyset$), then
$$
\Tilde H^i(\mathcal{C}^{\lpnt}_{\str_\Sigma C})
\cong
\Tilde H_{i-\dim C}(\lk_{\Delta(\Sigma)} \Fcc )
\cong
\begin{cases}
K  & \text{if } i=\dim C-1,\\
0  & \text{otherwise}.
\end{cases}
$$
The last isomorphism  follows from \cite[Proposition.4.7.8]{BJetal99}.

The supplement concerning the reduced simplicial homology  is obvious since we have that $\Delta(\str_\Sigma (0))=\Delta(\Sigma )$.
\end{proof}

We immediately get that the Cohen--Macaulay property of $K[\Sigma]$ only depends on the topology of
$|\Gamma_\Sigma|$.

\begin{corollary}
\label{examples}
Let $\Sigma$ be a $k$-dimensional rational pointed fan in $\RR^d$. The following conditions are equivalent:
\begin{enumerate}
\item[\rm{(a)}]
$K[\Sigma]$ is Cohen--Macaulay;
\item[\rm{(b)}]
For all $p\in |\Gamma_\Sigma|$ and all $i< k-1$ one has
$$
\tilde H_i(|\Gamma_\Sigma|)=H_i(|\Gamma_\Sigma|, |\Gamma_\Sigma| \setminus\{p\})=0.
$$
\end{enumerate}
Moreover, if the equivalent conditions are satisfied, then $\Sigma$ is Euler (and $K[\Sigma]$
is Gorenstein) if and only if for all $p\in |\Gamma_\Sigma|$
$$
\tilde H_{k-1}(|\Gamma_\Sigma|)=H_{k-1}(|\Gamma_\Sigma|, |\Gamma_\Sigma| \setminus\{p\})=K.
$$
In particular:
\begin{enumerate}
\item
If $|\Gamma_\Sigma|$ is homeomorphic to a $(k-1)$-dimensional sphere,
then $\Sigma$ is Euler and $K[\Sigma]$ is Gorenstein.
\item
If $|\Gamma_\Sigma|$ is homeomorphic to a $(k-1)$-dimensional ball,
then $K[\Sigma]$ is Cohen--Macaulay.
\end{enumerate}
\end{corollary}
\begin{proof}
(a)$\Leftrightarrow$(b): This is Corollary \ref{reisner} where one replaces  $\Tilde H^i(\mathcal{C}^{\lpnt}_{\str_\Sigma C})$
according to Lemma \ref{topolinterpreter}.
Using  Proposition \ref{needit}, (i)$\Leftrightarrow$(iv),
the remark concerning the Euler property follows.

(i):
If $|\Gamma_\Sigma|$ is homeomorphic to a $(k-1)$-dimensional sphere then
$$
\Tilde H_i(|\Gamma_\Sigma|;K)
=
H_i(|\Gamma_\Sigma|, |\Gamma_\Sigma| \setminus\{p\})
=
\begin{cases}
K & \text{if } i=k-1,\\
0 & \text{otherwise}.
\end{cases}
$$
for all  $p \in |\Gamma_\Sigma|$ (independent of the field $K$).

(ii): The second statement is shown analogously by observing that the local homology
$H_i(|\Gamma_\Sigma|, |\Gamma_\Sigma| \setminus \{p\})=0$ for $i\neq k-1$ if
$|\Gamma_\Sigma|$ is homeomorphic to a ball.
(We have that
$H_{k-1}(|\Gamma_\Sigma|, |\Gamma_\Sigma| \setminus \{p\})=0$ if $p$ is a boundary point.
Thus we only get the Cohen--Macaulay property.)
\end{proof}

Next we
give an application of Theorem \ref{main}.
For simplicial complexes the next theorem is due to Hochster.

\begin{theorem}
\label{hochster2}
Let $\Sigma$ be a rational pointed fan in $\RR^d$ such that $K[\Sigma]$ is a Cohen--Macaulay ring.
Assume that $X=|\Gamma_\Sigma|$
is a manifold with a non-empty boundary
$\partial X$.
Further
let $\Sigma'$ be the subfan of $\Sigma$
such that $\partial X=|\Gamma_{\Sigma'}|$.
Then the following conditions are equivalent:
\begin{enumerate}
\item $\omega_{K[\Sigma]}\iso\qq^\Sigma_{\Sigma'}$ as $\ZZ^d$-graded $K[\Sigma]$-modules;
\item $\Sigma'$  is Euler and $K[\Sigma']$ is Cohen--Macaulay.
\end{enumerate}
\end{theorem}
\begin{proof}
(i)$\Rightarrow$(ii):
Assume that $\Sigma$ is Euler. Then
$K[\Sigma]$ would be Gorenstein and
$K[\Sigma]=\omega_{K[\Sigma]}=\qq^\Sigma_{\Sigma'}$
would be a contradiction to the fact that
$\Sigma'$ is non-trivial.
Now (ii) was shown in (the proof of) Theorem \ref{main}.

(ii)$\Rightarrow$(i): In order to use Theorem \ref{main} we have to check that for all
$C\in\Sigma$
$$
 \Tilde H^{\dim\Sigma-1}(\mathcal{C}^{\lpnt}_{\str_\Sigma(C)}) =
\begin{cases}
0 & \text{ if }  C\in\Sigma',\\
K & \text{ if }  C\not\in\Sigma'.
\end{cases}
$$
We have to distinguish two cases:

(a) Let $\{0\} \neq C \in \Sigma$ and choose $p \in  e_C$.
We deduce that
$$
\Tilde H^{\dim\Sigma-1}(\mathcal{C}^{\lpnt}_{\str_\Sigma(C)})
\cong
H_{\dim\Sigma-1}(|\Gamma_\Sigma|, |\Gamma_\Sigma| \setminus\{p\})
=
\begin{cases}
0 & \text{ if }  p \in\partial X,\\
K & \text{ if }  p \not\in \partial X.
\end{cases}
$$
The first  isomorphism follows by Lemma \ref{topolinterpreter}, while the last equality
is implied by the definition of  a manifold with  boundary.
 Now $p \in\partial X$ if and only if $C \in \Sigma'$
and thus we are done in this case.

(b) Let $C=\{0\}  \in \Sigma'$. Then we have by Lemma \ref{topolinterpreter} that
$$
\Tilde H^{\dim\Sigma-1}(\mathcal{C}^{\lpnt}_{\str_\Sigma(C)})
\cong
H_{\dim\Sigma-1}(|\Gamma_\Sigma|).
$$
That $H_{\dim\Sigma-1}(|\Gamma_\Sigma|)=0$
was shown in the proof of \cite[Theorem 5.7.2]{BH}
applied to the simplicial complex $\Delta(\Sigma)$
with geometric realization $X=|\Gamma_\Sigma|$.
\end{proof}

%
%
%
\section{Shellability conditions and their algebraic properties}
\label{shellfan}

In \cite{IR07} the authors studied
the notion of a pure shellable fan and non-pure shellable fan respectively,
which imply that the corresponding
toric face ring is Cohen--Macaulay and sequentially Cohen--Macaulay respectively.
If the fan is the fan associated to a simplicial complex, then these notions
coincides with the well-known definitions of
pure shellable simplicial complexes and non-pure shellable simplicial complexes respectively.
For simplicial complexes there exist various equivalent definitions of shellability.

The goal of this section is to present several
shellability related notions for a rational pointed fan in $\RR^d$. Some of them
coincide in the ``simplicial case''. 

Let $\Sigma$ be a rational pointed fan in $\RR^d$.
In the following $\partial C$ denotes the boundary of a cone $C \in \Sigma$.
Recall that the fans $\fan(C)$ and $\fan(\partial C)$ are the set of faces of $C$ and $\partial C$ respectively.
Observe that they are subfans of $\Sigma$. At first we give the definition 
of shellability as considered in \cite{IR07} (which is motivated by the results of \cite{BJWA96} and \cite{BJWA97}).

\begin{definition}
A {\em shelling} of $\Sigma$
is a linear ordering $C_1,\dots, C_s$
of the facets of $\Sigma$ such that either $\dim \Sigma =0$,
or the following two conditions are satisfied:
\begin{enumerate}
\item
$\fan(\partial C_1)$ has a shelling.
\item
For $1<j\leq s$
there exists a shelling
$D_{1},\dots, D_{t_j}$ of the (pure) fan $\fan(\partial C_j)$
such that
$
\emptyset
\neq
\bigcup_{i=1}^{j-1} \fan(C_i) \cap \fan(C_j)
=
\bigcup_{l=1}^{r_j} \fan(D_l)$
for some $1\leq r_j \leq t_j$.
\end{enumerate}
$\Sigma$ is called {\em shellable} if it has a shelling.
\end{definition}

\begin{remark}
\
\begin{enumerate}
\item
Usually one considers only the case that $\Sigma$ is pure, i.e.\ all facets of $\Sigma$
have the same dimension.
But in this section we allow the fans also not to be pure
and all our shellings are non-pure shellings in the terminology of \cite{IR07} if not otherwise stated.
\item
One can weaken this definition in an obvious way by not asking
that the shelling of the fan in (ii) is the beginning of a shelling of
$\fan(\partial(C_j))$ and gets a weaker notion of shellability,
but we do not stress this point in this paper. Instead
we consider semishellings defined as below.
\end{enumerate}
\end{remark}

The next definition was inspired by  \cite{DK74}. (See also \cite{MCSH71}.)

\begin{definition}
A {\em semishelling} of $\Sigma$
is a linear ordering $C_1,\dots, C_s$
of the facets of $\Sigma$ such that either $\dim \Sigma=0$,
or
$|\Gamma_{\bigcup_{i=1}^{j-1} \fan(C_i) \cap \fan(C_j)}|$
is  homeomorphic either to a $(\dim C_j-2)$-dimensional ball or sphere for all $1< j \leq s$.
The fan $\Sigma$ is called {\em semishellable} if it has a semishelling.
\end{definition}

Of course we have:

\begin{corollary}\label{shellable_to_semishellable}
Let $\Sigma$ be a rational pointed fan in $\RR^d$.
If $\Sigma$ is shellable, then $\Sigma$ is semishellable.
\end{corollary}
\begin{proof}
Assume that $\dim \Sigma >0$ and
there is a linear ordering $C_1,\dots, C_s$
of the facets of $\Sigma$
such that
\begin{enumerate}
\item
$\fan(\partial C_1)$ has a shelling.
\item
For $1<j\leq s$
there exists a shelling
$D_{1},\dots, D_{t_j}$ of the (pure) fan $\fan(\partial C_j)$
such that
$
\emptyset
\neq
\bigcup_{i=1}^{j-1} \fan(C_i) \cap \fan(C_j)
=
\bigcup_{l=1}^{r_j} \fan(D_l)$
for some $1\leq r_j \leq t_j$.
\end{enumerate}
Now observe that the shellability of $\Sigma$ is equivalent to the fact
that then the regular cell complex $\Gamma_\Sigma$
is shellable in the sense of \cite[Section 4.7]{BJetal99}
with shelling order $e_{C_1},\dots, e_{C_s}$.
Note that
$|\Gamma_{\fan(\partial C_j)}|$ is a shellable $(\dim C_j-2)$-sphere for every $1\leq j \leq s$.
It follows now from \cite[Proposition 4.7.26 (i),(ii)]{BJetal99}
that
$|\Gamma_{\bigcup_{i=1}^{j-1} \fan(C_i) \cap \fan(C_j)}|$
is either homeomorphic to a $(\dim C_j-2)$-dimensional ball or sphere.
\end{proof}

Semishellability has the following nice algebraic consequence.
\begin{theorem}
\label{cm}
Let $\Sigma$ be a rational pointed fan in $\RR^d$.
If $\Sigma$ is pure semishellable, then $K[\Sigma]$ is Cohen--Macaulay
(independent of $\cha K$).
\end{theorem}
\begin{proof}
Using Corollary \ref{examples}
the proof of \cite[Theorem 3.2]{IR07} can be easily modified
to show that $K[\Sigma]$ is Cohen--Macaulay.
\end{proof}

Semishellability has also some consequences related to canonical modules.

\begin{proposition}
\label{weak}
Let $\Sigma$ be a rational pointed semishellable fan in $\RR^d$.
Then the toric face ring
$K[\bigcup_{i=1}^{j-1} \fan(C_i) \cap \fan(C_j)]$ is Cohen--Macaulay and
$
\omega_{\bigcup_{i=1}^{j-1} \fan(C_i) \cap \fan(C_j)}$
is isomorphic to a $\ZZ^d$-graded ideal of
$K[\bigcup_{i=1}^{j-1} \fan(C_i) \cap \fan(C_j)]$
for all $1< j \leq s$.
\end{proposition}
\begin{proof}
Assume that $\Sigma$ is semishellable with
semishelling order $C_1,\dots, C_s$.
Let $1< j \leq s$, $k=\dim C_j-1$
and consider the $k$-dimensional fan
$\Sigma'=\bigcup_{i=1}^{j-1} \fan(C_i) \cap \fan(C_j)$.
It follows from the definition of semishellability that
$|\Gamma_{\Sigma'}|$ is either homeomorphic to a $(k-1)$-dimensional ball or sphere.

If $|\Gamma_{\Sigma'}|$ is homeomorphic to a sphere, then $K[\Sigma']$ is Gorenstein
by Corollary \ref{examples}
and thus $\omega_{\Sigma'}=K[\Sigma']$.
Suppose that $|\Gamma_{\Sigma'}|$ is homeomorphic to a ball. Then $K[\Sigma']$ is Cohen--Macaulay.
There exists a subfan $\Sigma''$ of $\Sigma'$ such that
$|\Gamma_{\Sigma''}|$ is the boundary of the manifold with boundary $|\Gamma_{\Sigma'}|$ and hence homeomorphic to a sphere.
It follows from
Corollary \ref{examples} that $\Sigma''$ is Euler and that $K[\Sigma'']$ is Cohen--Macaulay.
Now it follows from Theorem \ref{hochster2} that
$\omega_{K[\Sigma']}$ is isomorphic to a $\ZZ^d$-graded ideal of
$K[\Sigma']$.
\end{proof}

\begin{remark} One can weaken the definition of a semishellable fan by asking only that
$$
\Tilde H_k(|\Gamma_{\bigcup_{i=1}^{j-1} \fan(C_i) \cap \fan(C_j)}|;K)
=
\begin{cases}
0 \text{ or } K & \text{if } 0\le k=\dim C_j-2,\\
0 & \text{otherwise}.
\end{cases}
$$
Theorem \ref{cm} and Proposition \ref{weak}  hold also in this context, and one may also give a converse for Proposition \ref{weak}.
\end{remark}

Next we consider a stronger property than shellability. For this we recall the following definition
which is due to Dress \cite{D}.

\begin{definition}
Let  $R$ be a Noetherian ring and $M$ be a finitely generated  $R$-module. A finite filtration
$$
0=M_0 \subset M_1 \subset \dots \subset M_r=M
$$
of submodules of $M$ is called a {\em prime filtration},
if for every $1\le i\le r$ there is an isomorphism $M_i/M_{i-1}\iso R/P_i$ for some prime ideal $P_i$ of $R$.
 A prime filtration
is called {\em clean}
if its set of corresponding prime ideals is equal to $\Min(\Supp M)$.
The module $M$ is called {\em clean}, if $M$ admits a clean filtration.
\end{definition}
Observe that there always exists a prime filtration, but clean filtrations may not exist.
For simplicial complexes,
cleanness is the algebraic counterpart of shellability as defined above.
Next we want to study  when $K[\Sigma]$ is clean for some rational pointed fan
$\Sigma$.
It is easy to see that if $K[\Sigma]$ is clean, then $\Sigma$ is shellable.
However, cleanness is a stronger property.

\begin{theorem}
\label{clean}
Let $\Sigma$ be a rational pointed fan in $\RR^d$.
Then $K[\Sigma]$ is clean if and only if the following conditions are satisfied:
\begin{enumerate}
\item $\Sigma$ is (non-pure) shellable with shelling $C_1,\ldots,C_s$.
\item For $1<j\le s$ there exists $\gamma_j\in C_j \cap \ZZ^d$ such that
\begin{enumerate}
\item $$\str_{\fan(C_j)}(\gamma_j)=\fan(C_j)\setminus\big(\bigcup_{i=1}^{j-1} \fan(C_i) \big).$$
\item $$\bigcup_{D\in \str_{\fan(C_j)}(\gamma_j)}\relint(D)\cap \ZZ^d=\gamma_j+C_j\cap\ZZ^d.$$
\end{enumerate}
\end{enumerate}
In particular, the cleanness property does not depend on $\cha K$.
\end{theorem}

\begin{proof} First assume that $K[\Sigma]$ is clean.
The minimal prime ideals of $K[\Sigma]$ are exactly the $\ZZ^d$-graded prime ideals $\pp_C$
for the maximal cones $C$ of $\Sigma$. Since $K[\Sigma]$ is reduced,
it follows from the last remarks in  \cite[Section 3]{D} that a clean prime filtration of $K[\Sigma]$
is necessarily of the form
$$
0=\bigcap_{i=1}^s\pp_{C_i} \subset \bigcap_{i=1}^{s-1} \pp_{C_i} \subset \dots \subset \pp_{C_1}\subset K[\Sigma]
$$
where $C_1,\ldots,C_s$ are the maximal cones of $\Sigma$, and there exist $\gamma_j\in C_j \cap \ZZ^d$
such that
$$
\bigcap_{i=1}^{j-1}\pp_{C_i}  \big/\bigcap_{i=1}^{j}\pp_{C_i}  \iso K[\fan(C_j)](-\gamma_j)\quad\text{ for }\quad 1< j\le s.
$$
Consider the fans $\Pi_j=\bigcup_{i=1}^{j-1} \big [\fan(C_i) \cap \fan(C_j)\big ]$ for $j=2,\dots,s$. Then
$$
\bigcap_{i=1}^{j-1}\pp_{C_i}  \big/\bigcap_{i=1}^{j}\pp_{C_i} =\qq_{\Pi_j}^{K[\fan(C_j)]}
\subseteq
K[\fan(C_j)]
\quad\text{ for }\quad 1< j\le s,
$$
so $\qq_{\Pi_j}^{K[\fan(C_j)]}\iso K[\fan(C_j)](-\gamma_j)=x^{\gamma_j}K[\fan(C_j)]$.
Let $D\in \fan(C_j)$. Then
$$
D\not\in \Pi_j\Longleftrightarrow
\relint(D)\cap \ZZ^d \subseteq \Supp(\qq_{\Pi_j}^{K[\fan(C_j)]},\ZZ^d).
$$
Assume $\gamma_j\in D$. Then $a+\gamma_j\in\relint(D)$
for all $a\in\relint(D)$.
It follows that
$(\relint(D)\cap \ZZ^d) \subseteq \Supp(K[\fan(C_j)](-\gamma_j),\ZZ^d)$.
Clearly $\gamma_j\not\in D$ implies $K[\fan(C_j)](-\gamma_j)\subset\pp_D^{K[\fan(C_j)]}$, so we have
$(D\cap \ZZ^d) \cap \Supp(K[\fan(C_j)](-\gamma_j),\ZZ^d)=\emptyset$.
We conclude that
$$D\in\str_{\fan(C_j)}(\gamma_j)\Longleftrightarrow
\relint(D)\cap \ZZ^d \subseteq \Supp(\qq_{\Pi_j}^{K[\fan(C_j)]},\ZZ^d)$$
which implies (ii) (a), and (ii) (b)
$$
\bigcup_{D\in \str_{\fan(C_j)}(\gamma_j)} \relint(D)\cap \ZZ^d=\Supp(\qq_{\Pi_j}^{K[\fan(C_j)]},\ZZ^d)=\gamma_j+C_j\cap\ZZ^d.
$$
It remains to show (i).
At first observe that $\fan(\partial C_j)$ has a shelling
because this is equivalent to the fact the boundary
of a cross-section polytope of $C_j$ has a shelling and this
is well-known by Bruggesser-Mani.
Using the notation introduced so far we have to prove
for $1<j\leq s$
that $\Pi_j=\bigcup_{l=1}^{r_j} \fan(D_l)$
where
$D_{1},\dots, D_{t_j}$ is a shelling of the fan
$\fan(\partial C_j)$ and $1\leq r_j \leq t_j$.

If $\gamma_j\in\relint(C_j)$, then
$\Pi_j=\fan(\partial C_j)$ and the assertion follows again
directly from Bruggesser--Mani.
Assume that $\gamma_j\not\in\relint(C_j)$.
Let $P_j$ be a cross-section polytope of $C_j$ such that
$\gamma_j \in P_j$.
It is possible to choose a point $x$ outside $P_j$
near $\gamma_j$ in general position with respect to $P_j$
(i.e.\ $x$ does not belong to an irredundant hyperplanes
defining $P_j$) such that the facets of $P_j$ which
correspond to the maximal cones in
$\str_{\fan(C_j)}(\gamma_j) \setminus \{C_j\}$
are exactly the visible facets from $x$.
(A facet is visible if for all points $y$ on that facet
the line segment between $y$ and $x$ does intersect $P_j$ only in $y$.)
Now it follows from  \cite[Theorem 8.12]{Z},
that there exists a (line) shelling of $P_j$
such that the facets  induced by the maximal cones in
$\str_{\fan(C_j)}(\gamma_j) \setminus \{C_j\}$
are the last ones.
This induces a shelling of $\partial C_j$
where
the maximal cones in
$\str_{\fan(C_j)}(\gamma_j) \setminus \{C_j\}$
are the last ones. Hence
the other ones which correspond to those cones
in $\Pi_j$ are the first ones.
Thus we have proved (i).

Now assume that (i) and (ii) holds. Consider the filtration of $K[\Sigma]$
$$
0=\bigcap_{i=1}^s\pp_{C_i} \subset \bigcap_{i=1}^{s-1} \pp_{C_i} \subset \dots \subset \pp_{C_1}\subset K[\Sigma]
$$
induced by the shelling of (i) and the fans $\Pi_j$ as above. Then
$$
\bigcap_{i=1}^{j-1}\pp_{C_i}  \big/\bigcap_{i=1}^{j}\pp_{C_i}  = \qq_{\Pi_j}^{K[\fan(C_j)]}
\subseteq
K[\fan(C_j)]\quad\text{ for }\quad 1< j\le s.
$$
It follows from (ii) that $\Supp(\qq_{\Pi_j}^{K[\fan(C_j)]},\ZZ^d)=\Supp(x^{\gamma_j}K[\fan(C_j)],\ZZ^d)$,
so the two ideals coincide. But
$x^{\gamma_j}K[\fan(C_j)]\iso K[\fan(C_j)](-\gamma_j)\iso \big(K[\Sigma]/\pp_{C_j}\big)(-\gamma_j),
$
hence $K[\Sigma]$ is clean.
\end{proof}

\begin{example}
It follows from Theorem \ref{clean} that
if $K[\Sigma]$ is clean, then $\Sigma$ is shellable.
As noted abover, it is known that for a simplicial complex $\Delta$
the converse is also true: If $\Delta$ is shellable,
then the Stanley--Reisner ring $K[\Delta]$ is clean.
(See Dress \cite{D}.)
For toric face rings in general the converse does not hold.
Indeed, e.g. consider the fan $\Sigma$ with facets $C_1,C_2 \subset \RR^2$
where
$$
C_1=\cn((0,1),(2,1))
\text{ and }
C_1=\cn((0,1),(-2,1)).
$$
Then it is easy to see that $\Sigma$ is shellable. The ring $K[\Sigma]$ is not clean
because it is not possible to find $\gamma_j$ as given in Theorem \ref{clean}.
\end{example}

\begin{proposition}
\label{conseq_clean}
Let $\Sigma$ be a rational pointed  fan in $\RR^d$ and let $K[\Sigma]$ be a clean ring. Let $C_1,\ldots,C_s$ be the
shelling of $\Sigma$ induced by a clean prime filtration of   $K[\Sigma]$.
By definition, for  $1<j\le s$
there exists a shelling
$D_{1},\dots, D_{t_j}$ of the (pure) fan $\fan(\partial C_j)$
such that
$
\Pi_j=
\bigcup_{i=1}^{j-1} [\fan(C_i) \cap \fan(C_j)]
=
\bigcup_{l=1}^{r_j} \fan(D_l)$
for some $1\leq r_j \leq t_j$.
Consider the fan $\Omega_j=\bigcup_{l=r_j+1}^{t_j} \fan(D_l)$.
Then $\Omega_j=\emptyset$ or $K[\Omega_j]$ is Gorenstein  for all $1<j\le s$ (independent of $\cha K$).
\end{proposition}
\begin{proof}Fix $j$ such that $1<j\le s$ and suppose that
$r_j\not= t_j$ (otherwise $\Omega_j=\emptyset$). In order to prove that $K[\Omega_j]$ is Gorenstein we check the assumptions of
\ref{gorensteinsecond}.
First of all $\Omega_j$ is pure shellable, so $K[\Omega_j]$ is Cohen--Macaulay by \ref{cm}.
Let $\gamma_j$ be as in Theorem \ref{clean} and consider $D\in \fan(C_j)$. 
As seen in the proof of  Theorem \ref{clean},
we have
$$
D\not\in \Pi_j\Longleftrightarrow D\in\str_{\fan(C_j)}(\gamma_j).
$$
In particular, $D_l\in\str_{\fan(C_j)}(\gamma_j)$ for all $r_j<l\le t_j$ We conclude that
$$\gamma_j\in \big [\bigcap_{C \in \Omega_j \text{ maximal}}C\big ]\cap\ZZ^d.$$
Observe that $\str_{\Omega_j}(\gamma_j)=\str_{\fan(C_j)}(\gamma_j)\setminus C_j$ since $\Omega_j\cup\Pi_j\cup C_j=\fan(C_j)$.
Also, for $a\in C_j\cap\ZZ^d$ we have that
$$
\gamma_j+a\in\relint(C_j)\cap \ZZ^d\Longleftrightarrow a\not\in\Supp(\Omega_j,\ZZ^d).
$$
By (ii)(b) in  Theorem \ref{clean} we have
\begin{align*}
\bigcup_{D\in \str_{\Omega_j}(\gamma_j)}\relint(D)\cap \ZZ^d
&=\Big[\bigcup_{D\in \str_{\fan(C_j)}(\gamma_j)}\relint(D)\setminus\relint(C_j)\Big]\cap \ZZ^d\\
&=(\gamma_j+C_j\cap\ZZ^d)\setminus (\relint(C_j)\cap \ZZ^d)\\
&=(\gamma_j+C_j\cap\ZZ^d)\setminus \Big\{\gamma_j+a:a\in (C_j\cap\ZZ^d)\setminus\Supp(\Omega_j,\ZZ^d)\Big\}\\
&=\gamma_j+\Supp(\Omega_j,\ZZ^d),
\end{align*}
so assumption (i) of \ref{gorensteinsecond} is satisfied.

It only remains to check assumption (ii) of \ref{gorensteinsecond}. We have to show that
for all cones $D \in \Omega_j$
we have
$$
\rho_{\Omega_j} (D)=
\begin{cases}
(-1)^{\dim \Omega_j - 1} &\text{ if } D \in \str_{\Omega_j}(\gamma_j) ;\\
0 &\text{ else. }
\end{cases}
$$
Since
$$
\rho_{\Omega_j} (D)=
\sum_{i=\dim D-1}^{\dim \Omega_j-1} (-1)^{i} \dim_K \Tilde H^{i}(\mathcal{C}^{\lpnt}_{\str_{\Omega_j}(D)}),
$$
and $K[\Omega_j]$ is Cohen--Macaulay, by Corollary \ref{reisner} it is enough to check that
$$
 \Tilde H^{\dim\Omega_j-1}(\mathcal{C}^{\lpnt}_{\str_{\Omega_j}(D)}) =
\begin{cases}
K & \text{ if }  D\in\str_{\Omega_j}(\gamma_j),\\
0 & \text{ else. }
\end{cases}
$$

Observe that $|\Gamma_{\Omega_j}|$ is homeomorphic to a
($\dim C_j-2$)-dimensional ball (as was shown in the proof of Corollary \ref{shellable_to_semishellable}).
There exists a subfan $\Lambda_j$ of $\Omega_j$ such that
$|\Gamma_{\Lambda_j}|$ is the boundary of the manifold with boundary $|\Gamma_{\Omega_j}|$ and hence homeomorphic to a sphere.
For $D\in\Omega_j$ we have that
$$
\Tilde H^{\dim\Omega_j-1}(\mathcal{C}^{\lpnt}_{\str_{\Omega_j}(D)})
=
\begin{cases}
K & \text{if }  D\not\in\Lambda_j,\\
0 & \text{if }  D\in\Lambda_j.
\end{cases}
$$
which one proves analogously to Corollary \ref{examples}
using the fact that one knows the local homology of $|\Gamma_{\Omega_j}|$.

Now consider $X=|\Gamma_{\fan(\partial C_j)}|$ as a topological space (it is a $(\dim C_j-1)$-sphere)
and $|\Gamma_{\Omega_j}|$ as a subspace in $X$. Then
$$
|\Gamma_{\Lambda_j}|=\partial |\Gamma_{\Omega_j}|=\overline{|\Gamma_{\Omega_j}|}\cap\overline{X\setminus|\Gamma_{\Omega_j}|}=|\Gamma_{\Omega_j}|\cap|\Gamma_{\Pi_j}|
$$
since $|\Gamma_{\Omega_j}|$ and $|\Gamma_{\Pi_j}|$ are ($\dim C_j-2$)-dimensional balls (so they are closed in $X$).

We deduce that $\Lambda_j=\Omega_j\cap\Pi_j$ and that for $D\in\Omega_j$ we have
$$
D\not\in \Lambda_j\Longleftrightarrow D\not\in \Pi_j\Longleftrightarrow D\in\str_{\Omega_j}(\gamma_j).
$$
This concludes the proof.
\end{proof}

\end{document}